\documentstyle[12pt]{article}
\normalbaselines
\begin{document}
\bibliographystyle{unsrt}
\vbox {\vspace{6mm}}
\begin{center}
{\large \bf ON SUMS OF POWERS OF ZEROS OF POLYNOMIALS \footnote {
Research talk given at the VIII Simposium sobre Polinomios Ortogonales y
Aplicaciones, Sevilla, Spain, September 22-27,1997 
} }\\
[9mm]
Wolfdieter L a n g \\ E-mail: wolfdieter.lang@physik.uni-karlsruhe.de\\
http://www-itp.physik.uni-karlsruhe.de/${\ \tilde{}}\ $wl\\ [3mm]
{\it Institut f\"ur Theoretische Physik \\ Universit\"at Karlsruhe \\
Kaiserstrasse 12, D-76128 Karlsruhe, Germany}\\[5mm]
\end{center}
\vspace{2mm}
\begin{abstract}
Due to Girard's (sometimes called Waring's) formula the sum of the 
$r-$th power of the zeros of every one variable polynomial of degree $N$, 
$P_{N}(x)$, can be given explicitly in terms of the coefficients of the monic 
${\tilde P}_{N}(x)$ polynomial. This formula is closely related to a known \par
\noindent $N-1$ variable generalization of Chebyshev's polynomials of the first kind, 
$T_{r}^{(N-1)}$. The generating function of these power sums 
(or moments) is known to involve the logarithmic derivative of the considered 
polynomial. This entails a simple formula for the Stieltjes transform of 
the distribution of zeros. Perron-Stieltjes inversion can be used to find this
distribution, {\it e.g.} for $N\to \infty$.\par
Classical orthogonal polynomials are taken as examples. The results for
ordinary Chebyshev $T_{N}(x)$  and $U_{N}(x)$ polynomials are presented in 
detail. This will correct a statement about power sums of zeros of 
Chebyshev's $T-$polynomials found in the literature. For the various cases 
(Jacobi, Laguerre, Hermite) these moment generating functions provide 
solutions to certain Riccati equations. \par
\vskip 1cm  \par \noindent
\end{abstract}
\section{Introduction}
\hskip 1cm 
Sums of powers of the zeros of every one variable polynomial $P_{N}(x)$ of 
degree $N$ can be given in terms of its coefficients as $(N-1)$-fold nested sums.
If the polynomial at hand has definite parity, {\it i.e.} if $P_{N}(-x)\ =
(-1)^{N}P_{N}(x)\ $, only $\lfloor {{N}\over{2}}-1\rfloor $ nested sums are
needed. These facts follow from {\it Girard'}s formula \cite {McM 60}: The 
considered $r^{th}$ power sums are one part symmetric functions (polynomials) 
of weight $r$ in the zeros, and due
to the fundamental theorem of symmetric functions, {\it e.g.}, 
\cite{Kr 86} they can be uniquely expressed in terms of the elementary
symmetric functions $\{\sigma_{i}\}_{1}^{N}\ $ of the zeros, {i.e.} 
(modulo signs) in terms of the coefficients of the corresponding
monic polynomial ${\tilde P}_{N}(x)$. {\it Girard}'s formula which solves this
problem is usually, {\it e.g.}, \cite {McM 60} given in a 'semi-explicit' form 
where still partitions have to be performed. A simplification of this formula
has been given in \cite {LiWe 72},
\cite {Li 75}, \cite{LiMuTu 92}. Here we shall present the explicit formula. The latter authors 
defined a certain $(N-1)$ variable generalization of {\it Chebyshev}'s polynomials
of the first kind, $T_{r}^{(N-1)}(\sigma_{1},...,\sigma_{N})$, 
by the {\it Girard-Waring} expression after
putting $\sigma_{N}$, the product of the zeros, to one.\par\smallskip\noindent
The zeros of characteristic polynomials $P_{N}(x)$ determine the eigenvalue 
spectrum of (finite-dimen- \par\noindent
sional) operators. Therefore information on these 
zeros is welcome. In most cases these zeros are not known, and numerical or
perturbative methods are employed. Often, the limit case $N\to \infty$ is of
interest and the limit density of zeros 
$\rho(x)\ :=\lim_{N\to \infty} \rho(N,x)\ $ is asked for. Here the power sums
of zeros become important because, for given $P_{N}(x)$ with only real zeros 
(hence for all orthogonal polynomials with positive moment functional), the 
(ordinary) generating function of all (positive integer) power sums, 
$G(N,z)$, is related in a simple way to the {\it Stieltjes}-transform of the 
discrete distribution of zeros $\rho(N,x)$ 
(the counting measure of the zeros which lives on the real line). This 
function $G(N,z)$, which generates the generalized {\it Chebyshev} polynomials of the 
first kind, is essentially given by the logarithmic derivative of the polynomial
$P_{N}(x)$ under consideration. If it can be considered as a 
function of real 
$N$, and the limit $N\to\infty$ exists, then one can obtain the limit distribution
of zeros of $P_{N}(x)$ by {\it Perron-Stieltjes} inversion \cite{ Wi 29},\ 
\cite {Ak 65}, \cite {Ch 78},\ \cite{As 84}. \par \smallskip
In this report which is based on the more detailed work \cite {La 97} (which is 
referred to for proofs and intermediate steps) we demonstrate the above 
mentioned facts for the classical orthogonal polynomials.
The generating functions $G(N,z)$ are computed, and they are seen to satisfy
certain {\it Riccati}-equations, or their equivalent linear 
second order differential equations. For {\it Chebyshev}'s polynomials of both 
kinds simple expressions for this generating functions of the moments of the 
dicrete counting measure of the zeros result. It becomes 
transparent why in these cases the $1/N$ expansion for the power sums or 
moments becomes unreliable. Finally the well-known limit distribution of zeros 
for the {\it Jacobi}-polynomials and (appropriately scaled) generalized 
{\it Laguerre-} and {\it Hermite-}polynomials are found from the 
{\it Stieltjes-}transform {\it via} {\it Perron-Stieltjes}-inversion. 
\par \smallskip
As explained above, this method to find limit distributions of zeros works 
for all polynomials with real zeros whose logarithmic 
derivative can be considered as function of the degree number $N$ with existing 
limit $N\to \infty$ (possibly rescaled).
\par
\bigskip\noindent
\section{Girard-Waring formula for the power sums of zeros of polynomials}
\par
A monic polynomial of degree $N$ in one real variable $x$ is written as  
\begin{equation}
\tilde{P}_{N}(x)\ =\sum_{k=0}^{N}\ (-1)^{k}\ \sigma_{k}(N)\ x^{N-k}\ \ \ \ = \ \ 
\ \prod_{i=1}^{N}\ (x-x_{i}^{(N)}) . 
\end{equation}
The fundamental theorem of algebra has been used. The second equation 
defines the $k^{th}$ elementary symmetric polynomials 
$\sigma_{k}(N)\equiv \sigma_{k}\ $
of the $N$ zeros $\{x_{i}^{(N)}\}_{1}^{N}$ of $P_{N}(x)$ .\  
$\sigma_{0}(N)\equiv 1\ . $
\par\smallskip\noindent
The normalized sum of the $r^{th}$ power of the zeros is
\begin{eqnarray}
m_{r}(N)&\ :=&{{1}\over{N}}\ \sum_{i=1}^{N}\ (x_{i}^{(N)})^{r} \\
&=& \int_{-\infty}^{+\infty} x^{r}\ d\rho_{N}(x)  \\
d\rho_{N}(x)&\ :=&{{1}\over{N}}\sum_{i=1}^{N} \delta(x-x_{i}^{(N)})\ dx .
\end{eqnarray}
\noindent The discrete counting measure for the zeros, $\rho_{N}(x)$, is 
used only for polynomials with only real zeros, in order to have this  
measure normalized to $1$. Under this assumption the power sums are the moments
of this real discrete measure. (Sometimes we shall call the power sums of 
zeros of {\it general} polynomials also moments.)\par
For an arbitrary polynomial {\it Girard}'s  formula, when written explicitly,
looks  for $r\in {\bf N}\ \ , \ \ N= 2,3,...$ as follows:
\begin{eqnarray}
m_{r}(N)\ = &{{1}\over{N}}\ 
\sigma^{r}_{1}\ \Bigl{(} \prod_{k=1}^{N-1} \bigl{(} \sum_{i_{k}=0}
^{\lfloor b_{k} \rfloor}\ {{(-1)^{(N-k)i_{k}}}\over{(i_{k})!}}\bigr{)}\Bigr{)}\ 
 <N,r,\{i_{j}\}_{1}^{N-1}>\ \prod_{p=1}^{N-1} \bigl{(}{{\sigma_{p+1}}\over
{\sigma_{1}^{p+1}}}\bigr{)}^{i_{N-p}}\ \ &  \\
&&\nonumber \\
&=: \ t_{r}^{(N)}(\sigma_{1},...,\sigma_{N}) \ \ , & 
\end{eqnarray}
\noindent 
where the $N-1$ sums are ordered from the left to the right with increasing $k$
value, and the upper boundaries of these sums are the greatest integer of
$b_{1}:=r/N$  and
\begin{equation}
b_{k}\ := {{r\-\sum_{j=1}^{k-1}\ (N+1-j)\ i_{j}}\over{N+1-k}}\ \ \ , \ \ 
k=2,...,N-1 \ . 
\end{equation}
The combinatorial factor which appears in eq. $(5)$ is
\begin{equation}
<N,r,\{i_{j}\}_{1}^{N-1}>\ := 
{{r\ (r-1-\sum_{j=1}^{N-1}\ (N-j)i_{j})!}
\over{(r-\sum_{j=1}^{N-1}\ (N+1-j)\ i_{j})!}}\ . 
\end{equation}
\noindent For $r=0$ one has $m_{0}(N)\equiv 1$. 
$t_{r}^{(N)}(\sigma_{1},...,\sigma_{N})$, which is a polynomial, has been 
introduced because it is related to an $N-1$ variable generalization of
{\it Chebyshev}'s polynomial of the first kind $T_{r}(x)$. This connection is
given after the rescaling 
\begin{equation}
T_{r}^{(N-1)}(s_{1},...,s_{N-1})\ = t_{r}^{(N)}(\sigma_{1},...,\sigma_{N})/
(\sigma_{N})^{r/N}\ \ ,
\end{equation}
\noindent with
\begin{equation}
s_{k}\ := \sigma_{k}/(\sigma_{N})^{k/N}\ \ ,\ \ \ k=1,...,N-1\ .
\end{equation}
The ordinary {\it Chebyshev} polynomials are now $T_{r}(x)\ = T_{r}^{(1)}(2x)$.
\par\smallskip\noindent
{\bf Corollary: Moments for $\bf P_{N}(-x)=(-1)^{N}\ P_{N}(x)\ $ polynomials}
\par\smallskip\noindent
\begin{eqnarray} 
&i)\ \ \  m_{0}(N)\equiv 1\ \ ,\ \ N\in {\bf N_{0}}\ \hskip 2cm 
ii)\   m_{2l+1}(N)\equiv 0\ \ ,\ \ l\in {\bf N_{0}}\ ,\  N\in {\bf N}\ , & \\
&iii)\  m_{2l}(N)\ = & \nonumber \\
  &{{2}\over{N}}\ (-\sigma_{2})^{l}\ \Bigl{(}\prod_{k=1}
 ^{\lfloor {{N}\over{2}}\rfloor-1}\ \bigl{(}
 \sum_{i_{k}=0}^{\lfloor B_{k}\rfloor}\ 
 (-1)^{(\lfloor {{N}\over{2}} \rfloor-k)i_{k}}{{1}\over{(i_{k})!}}\bigr{)} 
 \Bigr{)}\ <\lfloor {{N}\over{2}}\rfloor,l,
 \{i_{j}\}_{1}^{\lfloor {{N}\over{2}}\rfloor-1}>\ 
\prod_{p=2}^{\lfloor {{N}\over{2}}\rfloor} \bigl{(}
{{\sigma_{2p}}\over{\sigma_{2}^{p}}}\bigr{)}^{i_{\lfloor {{N}\over{2}}\rfloor
-p+1}}\ \ & \\
&\ ={{2}\over{N}} \lfloor {{N}\over{2}}\rfloor (-1)^{l}\ 
t_{l}^{(\lfloor {{N}\over{2}}\rfloor)}(\sigma_{2},\sigma_{4},...,
\sigma_{2 \lfloor {{N}\over{2}}\rfloor}) &\\
&\ = {{2}\over{N}} 
\lfloor {{N}\over{2}}\rfloor (-1)^{l}\ 
\bigl{(}\sigma_{2 \lfloor {{N}\over{2}}\rfloor}
\bigr{)}^{l/ \lfloor {{N}\over{2}}\rfloor}\ 
T_{l}^{(\lfloor {{N}\over{2}}\rfloor -1)}(e_{1},...,
e_{\lfloor {{N}\over{2}}\rfloor -1})\ \ \ & \\
&&\nonumber  \\
&{\rm for}\  l\in {\bf N_{0}}\ , \ \ N=4,5,..., & \nonumber  
\end{eqnarray}
\noindent where 
$B_{k}$, resp. $<\lfloor{{N}\over{2}}\rfloor,l,\{i_{j}\}_{1}^{ \lfloor {{N}\over{2}}
\rfloor  -1}>\ $, resp. $t_{l}^{(\lfloor {{N}\over{2}}\rfloor)}$,
is obtained from $b_{k}$ of $(7)$, resp.\par \noindent 
$<N,r,\{i_{j}\}_{1}^{N-1}>\ $ of $(8)$, resp. $t_{r}^{(N)}$ of $(6)$, 
by replacing $N\to \lfloor {{N}\over{2}}\rfloor $ 
and $r\to l$. The variables of $T_{l}^{(\lfloor {{N}\over{2}}\rfloor -1)}$, 
which
is also obtained by the same replacement from $T_{r}^{(N-1)}$ given in $(9)$,
are
\begin{equation}
e_{k}\equiv e_{k}(N)\ =\sigma_{2k}/
\bigl{(}\sigma_{2\lfloor {{N}\over{2}}\rfloor}
\bigr{)}^{k/\lfloor {{N}\over{2}}\rfloor}\ \ \ ,\ \ \ 
k=1,...,\lfloor {{N}\over{2}}\rfloor -1 \ \ . 
\end{equation}
In addition: $m_{2l}(1)\ =0  \ \ \ , \ \ \ m_{2l}(2)\ = (-\sigma_{2})^{l}\ \ \
\ \ ,\ \ \ m_{2l}(3)\ = {{2}\over{3}}(-\sigma_{2})^{l} $.
\par\smallskip\noindent
\noindent{\bf Example:}\par\noindent
Consider {\it Chebyshev}'s polynomials of the first kind, $T_{N}(x)$, which are 
classical orthogonal polynomials of the {\it Jacobi} type \cite{Sz 59}, \cite{Ch 78}. 
The zeros are known, and the non-vanishing moments are 
$m_{2l}^{(T)}(N)\ ={{2}\over{N}}
\sum_{k=1}^{\lfloor {{N}\over{2}}\rfloor}
(cos \bigl{(}{{2k-1}\over{2N}}\pi\bigr{)})^{2l}\ $. The {\it Corollary} 
gives an alternative expression for 
this sum. \par\noindent
\begin{eqnarray}
&{\rm For\ } l\in {\bf N}\ ,\ M\in {\bf N}\ :\ 
m_{2l}^{(T)}(2M)\ =T_{l}^{(M-1)}(e_{1}(2M),...,
e_{M-1}(2M))/2^{(2-{{1}\over{M}})l} \ \ , & \\
&&\nonumber \\
&\hskip 1cm {\rm with}\ \ e_{k}(2M)\equiv {{2M}\over{k}}\ 2^{-{{k}\over{M}}} 
{{2M-k-1}\choose{k-1}}\ \ ,\ \ \ k=1,...,M-1, & \\
&m_{2l}^{(T)}(2M+1)\ ={{2M}\over{2M+1}}\ \bigl{(}{{(2M+1)^{{{1}\over{M}}}}\over
{4}}\bigr{)}^{l}\ T_{l}^{(M-1)}(e_{1}(2M+1),...,e_{M-1}(2M+1)) , & \\
& & \nonumber \\
&\hskip 1cm {\rm with}\ \ e_{k}(2M+1)\equiv {{2M+1}\over{k}}\ 
{{1}\over{(2M+1)^{{{k}\over{M}}}}}\ {{2M-k}\choose{k-1}} ,\  k=1,...,M-1.\ &
\end{eqnarray}
\noindent and $m_{0}^{(T)}(N)\ =1\ \ \ ,\ \ \ m_{2l}^{(T)}(1)\ =0\ $.\par
\noindent Some instances for $l\in {\bf N}$ are:
\begin{eqnarray}
&m_{2l}^{(T)}(2)\ = ({{1}\over{2}})^{l}\ \ ,\ \ m_{2l}^{(T)}(4)\ = T_{l}(\sqrt{2})/
(2\sqrt{2})^{l}\ \ , & \nonumber \\
&m_{2l}^{(T)}(6)\ =T_{l}^{(2)}(3\cdot 2^{{{2}\over{3}}},
9\cdot 2^{-{{2}\over{3}}})/2^{{{5}\over{3}}l}\ \ ,   
& \\
&m_{2l} ^{(T)}(3)\ =3^{l-1}/2_{2l-1}\ \ \ ,\ \ \ m_{2l}^{(T)}(5)\ =
(\sqrt{5})^{l-2}\ T_{l}(\sqrt{5}/2)/2^{2(l-1)}\ \ , \ \ & \nonumber \\
&m_{2l}^{(T)}(7)\ = {{6}\over{7}}\ \bigl{(}{{7^{{{1}\over{3}}}}\over{4}}
\bigr{)}^{l}\ T_{l}^{(2)}(7^{{{2}\over{3}}}, 2\cdot7^{{{1}\over{3}}})\ \ .
& 
\end{eqnarray}
From this one finds for $N$ up to  $7$:
\begin{equation}
m_{2}^{(T)}(N)\ =\cases{0&N=1\cr{{1}\over{2}}&N=2,...}\ \  ,\ \ 
m_{4}^{(T)}(N)\ =\cases{0&N=1\cr1/4&N=2\cr3/8&N=3,...\cr}\ \ \ ,
\ \  m_{6}^{(T)}(N)\ =\cases{0&N=1\cr1/8& N=2\cr9/32 &N=3\cr
5/16& N=4,...\cr}\ \ .
\end{equation}
\noindent These values are in conflict with a statement found in \cite {Ca 80}.
There it is claimed that the moments (22) are given by the highest values, 
independently of $N$, which is false. Only the values for all $N\geq {{r}\over{2}}+1$
are correctly given in \cite{Ca 80}.  \par\noindent
The example of {\it Chebyshev}'s polynomials of the second type is also treated
in \cite{La 97}.\par\bigskip\noindent
\section{Moment generating functions}\par\smallskip
The generating function of $p_{r}(N)\ $, the sum of the $r^{th}$ power of $N$ variables 
$\{x_{i}\}_{1}^{N}\ $, is given by the logarithmic derivative of the generating
function of the elementary symmetric functions $\{\sigma_{k}\}\ $ of these
variables, see {\it e.g.} \cite{Kr 86}. With 
\begin{equation}
E(z)\equiv E_{N}(z):=\prod_{i=1}^{N}\ (1-x^{(N)}_{i}\ z)\ = 
\sum _{r=0}^{N}\ (-1)^{r}\sigma_{r}\ z^{r}\
\ \  \ , \ \sigma_{0}=1\ \ , 
\end{equation}
which is also given by $E(z)\ = z^{N}{\tilde P}_{N}(1/z)\ $, one finds from 
\begin{equation}
(ln\ E(z))^{\prime}\ = E^{\prime}(z)/E(z)\ =-{{1}\over{z}}\sum_{r=1}^{\infty}
p_{r}(N)\ z^{r}\ \ 
\end{equation}
and {\it Girard's} fomula $(5),\  (6)$
\begin{eqnarray}
G(N,z)\ := \sum_{r=0}^{\infty}m_{r}(N)\ z^{r}\ &=\  
\sum_{r=0}^{\infty}\ t_{r}^{(N)}(\sigma_{1},...,\sigma_{N})\ z^{r}\ 
= \left.{{1}\over{N}} {{x {{d{\ }}\over{dx}} P_{N}(x)}\over{P_{N}(x)}}
\right\vert_{x=1/z} \ = \nonumber \\
&=\ {{\sum_{r=0}^{N-1}(-1)^{r}(1-{{r}\over{N}})\ 
\sigma_{r}(N)\ z^{r} }\over{ \sum_{r=0}^{N}(-1)^{r}\sigma_{r}(N)\ z^{r}
}} .
\end{eqnarray}
\noindent As an application we computed the moment generating functions $G(N,z)$ for
the classical orthogonal polynomials: {\it Jacobi} $P_{N}^{(\alpha,\beta)}(x)\ $,
generalized {\it Laguerre} $L_{N}^{(\alpha)}(x)\ $  and {\it Hermite} 
$ H_{N}(x)\ $ . The derivative rules for these polynomials ({\it cf.} \cite{Ch 78})
have been used. We quote the results for {\it Chebyshev}'s polynomials of both
kinds separately because they can be written in terms of an elementary 
function. For details see \cite{La 97}. The hypergeometric function $_{2}F_{1}
(a,b,c;x)\ $ and its confluent version $F(a,c;x)\equiv\  _{1}F_{1}(a,c;x)\ $ 
appear in {\it table 1}.\par\smallskip\noindent
The $(N+1)-$term recursion relation for the $t_{r}^{(N)}-$polynomials which 
can be derived from eq.(25), is:
\begin{equation}
{\rm for}\ n\ =\ N,N+1,...,:\ \ \ \ \sum_{j=0}^{N}(-1)^{j}\sigma_{j}\ 
t^{(N)}_{n-j}\ = \ 0 \ \ \ , \ \ \ \sigma_{0} =1\ , 
\end{equation}
where the $N$ input quantities $\{t_{j}^{(N)}\}_{0}^{N-1}\ $ are either 
computed from {\it Girard'}s formula or from the generating function $G(N,z)$.
\par
We also comment on the generating functions $G_{<}(N,z)\ $ for the normalized 
sums of negative powers of the zeros: $m_{-r}(N)\ :={{1}\over{N}} \sum_{i=1}^{N}\ 
(x_{i}^{(N)})^{-r}\ $, $r\in {\bf N}$. Here we restrict the discussion to polynomials 
which have no vanishing zeros. With $E_{N}^{<}(z)\ :=\ \prod_{i=1}
^{N}(1-z/x_{i}^{(N)})\ =$ \ $(-1)^{N} \tilde{P}_{N}(z)/\sigma_{N}\ $ one finds
\begin{equation}
G_{<}(N,z)\ := \ \sum_{r=1}^{\infty}\ m_{-r}(N)\ z^{r}\ = \ -{{z}\over{N}}
(ln\ E_{N}^{<}(z))^{\prime}\ 
= \ - {{z}\over{N}}{{\tilde{P}^{\prime}_{N}(z)}\over{\tilde{P}_{N}(z)}}\ =
\ -\ G(N,1/z)\ .
\end{equation}
Beware of expanding the identity $G_{<}(N,z) \ +\ G(N,1/z)\ \equiv \ 0$ 
in $z$ {\it and} $1/z\ $ (and concluding falsely that all power sums have to vanish). 
The sums of negative powers of the $P_{N}(x)$ zeros can then be found by 
differentiation. The connection between 
{\it Fa\`a di Bruno}'s formula and the {\it Bell-}polynomials 
\cite {Ri 58} and the 
identification of the generalized {\it Chebyshev} polynomials $t^{(N)}_{r}$ as 
special {\it Bell-}polynomials leads, with the abbreviation 
$g(x)\ \equiv \tilde{P}_{N}(x) \ $, for $m\in {\bf N}$ first to
\begin{eqnarray}
&{{d^{m}\ }\over{dx^{m}}}\ ln\ g(x)\ =  \nonumber \\
 &\ -(m-1)!\ N\ t^{(N)}\bigl{(}- g^{(1)}(x)/
(1!g(x)),\ +g^{(2)}(x)/(2!g(x)),\ ..., (-1)^{m}g^{(m)}(x)/(m!g(x))\bigr{)}\ ,
\end{eqnarray}
with $g^{(k)}\equiv {{d^{k}\ }\over{dx^{k}}} g(x)\ $. After putting $x=0$, 
using $(-1)^{k} g^{(k)}(0)/(k!g(0))\ = \ \sigma_{N-k}/\sigma_{N}\ $, one 
finally finds
\begin{equation}
m_{-r}(N)=\left\{ \begin{array}{ll}
t_{r}^{(N)}\bigl{(}\sigma_{N-1}/\sigma_{N},\ \sigma_{N-2}/ 
\sigma_{N},\ ...,\  \sigma_{N-r}/\sigma_{N},0,...,0\bigr{)} & \mbox{for}\  
\ 0\ \leq\  r\ \leq N  \\
&\\ 
t_{r}^{(N)}\bigl{(}\sigma_{N-1}/\sigma_{N},\ \sigma_{N-2}/ \sigma_{N},\ ...,\ 
(\sigma_{0}=1)/\sigma_{N}\bigr{)} &  \mbox{for}\  r\ \geq\  N \ , 
\end{array}\right. 
\end{equation}
where in the upper alternative $N-r$ zeros have to be inserted. \par \smallskip
\noindent 
\begin{center}
{\bf Table 1\ :  Generating functions for power sums of zeros of 
classical orthogonal polynomials}
\end{center}
\smallskip \noindent
\begin{center}
\begin{tabular}{|c|c|}\hline
& \\
$\bf P_{N}(x)$ & $\bf G(N,z)$ \\ 
&  \\ \hline\hline
& \\
$P^{(\alpha,\beta)}_{N}$ & $G^{(\alpha,\beta)}(N,z)\ =\ {{N+\alpha+\beta+1}
\over{2(\alpha+1)}}{{1}\over{z}}{{_{2}F_{1}(-N+1,\ N+\alpha+\beta+2,
\ \alpha+2;\ (z-1)/2z)}\over{_{2}F_{1}(-N,\ N+\alpha+\beta+1,\ 
\alpha+1;\ (z-1)/2z)}}$  \\ 
&\\ \cline{2-2}
&\\
$T_{N}$  &$G^{(T)}(N,z)\ =\ {{1}\over{\sqrt{1-z^2}}}tanh\bigl{(}N\ 
ln{{z}\over{1-\sqrt{1-z^2}}}\bigr{)}$\\
& \\ \cline{2-2}
&\\
$U_{N}$  & $G^{(U)}(N,z)\ =\ {{1}\over{\sqrt{1-z^2}}}\Bigl{\{}(1+{{1}\over{N}})\ tanh\Bigl{(}
(N+1)\ ln{{z}\over{1-\sqrt{1-z^2}}}\Bigr{)}-{{1}\over{N}}{{1}\over
{\sqrt{1-z^2}}}\Bigr{\}}$\\
& \\ \hline
& \\
$L^{(\alpha)}_{N}$  & $G^{(\alpha)}(N,z)\ =\
1-{{F(-N+1,\alpha+1;1/z)}\over{F(-N,\alpha+1;1/z)}}$ \\
&\\ \hline
& \\
$H_{N}$  & $ G^{(H)}(N,z)\ =\ {{1}\over{z}}{{F(-{{N-1}\over{2}},\ {{1}\over{2}}
;\ 1/z^2)/\Gamma({{2-N}\over{2}})\ - {{2}\over{z}}\ F({{2-N}\over{2}},\ 
{{3}\over{2}};\ 1/z^2)/\Gamma({{1-N}\over{2}})}\over{F(-{{N}\over{2}},
{{1}\over{2}};1/z^2)/\Gamma({{1-N}\over{2}})\ - 
{{2}\over{z}}\ F({{1-N}\over{2}},\  {{3}\over{2}};\ 1/z^2)/\Gamma(-{{N}\over{2}})
}}$ \\
&\\ \hline
\end{tabular}
\end{center}
\bigskip\noindent
One can verify ({\it cf.} \cite {La 97}) that for orthogonal polynomials the 
generating 
functions $G(N,z)$ are special solutions of certain {\it Riccati}-equations
which we list in the following {\it table 2}.
\vfill\eject\noindent
\begin{center}
{\bf Table 2\ : Riccati-equations satisfied by moment generating functions of 
classical orthogonal polynomials }
\end{center}
\bigskip \noindent
\begin{center}
\def\Paeq#1#2{{{\partial {#1}}\over{\partial {#2}}}}
\begin{tabular}{|c|c|}\hline
& \\ 
$\bf P_{N}(x)$& Riccati-eq. for $\bf G(N,z)$ for every $\bf N$ \\
& \\ \hline\hline
& \\
$P^{(\alpha,\beta)}_{N}$ & $ {{(1-z^{2})z}\over{N}}\Paeq{}{z}
G^{(\alpha,\beta)}(N,z) \ =
(z^{2}+(\alpha -\beta)z+\alpha +\beta+ 1){{1}\over{N}}\  
G^{(\alpha,\beta)}(N,z)$ \\
& \\
&$\ +(1-z^{2})\ \bigl{(} G^{(\alpha,\beta)}(N,z)\bigr{)}^{2}\ -
(1+{{1+\alpha +\beta}\over{N}})$   \\ 
& \\ \cline{2-2}
& \\
$T_{N}$ & ${{(1-z^2)z}\over{N}}\Paeq{}{z}G^{(T)}(N,z)\ = {{z^{2}}\over{N}}\  G^{(T)}(N,z)
\ +(1-z^{2})\bigl{(} G^{(T)}(N,z)\bigr{)}^{2}\ -1\ \ $ \\ 
& \\ \cline{2-2}
& \\
$U_{N}$ & $ {{(1-z^2)z}\over{N}}\Paeq{}{z}G^{(U)}(N,z)\ = 
{{2+z^{2}}\over{N}}\  G^{(U)}(N,z) \ +(1-z^{2})\bigl{(} G^{(U)}(N,z)
\bigr{)}^{2}\ -{{N+2}\over{N}} $ \\ 
 & \\ \hline
 & \\
$L^{(\alpha)}_{N}$ & $z^{2}\ \Paeq{}{z}\  G^{(\alpha)}(N,z)\ = Nz\bigl{(}
G^{(\alpha)}(N,z)\bigr{)}^{2}\ +(\alpha z-1)\ G^{(\alpha)}(N,z)\ +1$ \\ 
& \\ \cline{2-2}   
& \\
$Ls_{N}^{(\alpha)}$ &  $z^{2}{{1}\over{N}}\ 
\Paeq{}{z}\  Q^{(\alpha)}(N,z)\ = z\bigl{(}
Q^{(\alpha)}(N,z)\bigr{)}^{2}\ +(\alpha z/N-1)\ Q^{(\alpha)}(N,z)\ +1$ \\ 
& \\ \hline
& \\
$H_{N}$ &  $ {{z^{3}}\over{2}}\ \Paeq{}{z}\  G^{(H)}(N,z)\ = {{z^{2}}\over{2}}N\bigl{(}
 G^{(H)}(N,z)\bigr{)}^{2}\ -(1+{{z^{2}}\over{2}}) G^{(H)}(N,z)\ +1 $ \\ 
 & \\ \cline{2-2}
 & \\
$Hs_{N}$ &  $ {{z^{3}}\over{2N}}\ \Paeq{}{z}\  Q^{(H)}(N,z)
\ = {{z^{2}}\over{2}}\bigl{(}
Q^{(H)}(N,z)\bigr{)}^{2}\ -(1+{{z^{2}}\over{2N}}) Q^{(H)}(N,z)\ +1 $\\ 
 & \\ \hline  
\end{tabular}
\end{center}
\par\bigskip\noindent
We have also listed the {\it Riccati}-equations satisfied by the moment 
generating functions of scaled generalized {\it Laguerre-} and 
{\it Hermite-}polynomials.
The definitions are:
\begin{eqnarray}
Ls^{(\alpha)}_{N}(x)\ :=&L_{N}^{(\alpha)}(Nx) \hskip 1cm ,\hskip 1cm  
Q^{(\alpha)}(N,z)\ :=G^{(\alpha)}(N,z/N)\ & ,\\
Hs_{N}(x)\ :=&H_{N}(\sqrt{N}x) \hskip 1cm ,\hskip 1cm
Q^{(H)}(N,z)\ :=G^{(H)}(N,z/\sqrt{N})\ \ . 
\end{eqnarray}
\def\Paeq#1#2{{{\partial {#1}}\over{\partial {#2}}}}
\noindent The associated second order linear differential eqs. are listed
in {\it table 3}. $ G(N,z)\ $ , \par\noindent
${\rm resp.\ } Q^{(\alpha)}(N,z):=
G^{(\alpha)}(N,z/N)$ and $Q^{(H)}(N,z):=G^{(H)}(N,z/\sqrt{N})$, are replaced
by \par\noindent
$-{{z}\over{N}}\ \Paeq{}{z} \bigl{(}
ln\ {\cal U}(N,z)\bigr{)}$, with correspondingly labelled ${\cal U}(N,z)$.    
\vfill\eject\noindent
\begin{center}
{\bf Table\ 3: Second order linear differential equation satisfied by 
$\bf {\cal U}(N,z)$ for given $\bf N$ 
for (scaled) classical orthogonal polynomials}
\end{center}
\begin{center}
\begin{tabular}{|c|c|}\hline
&\\
$\bf P_{N}(x)$& $\bf 2^{\rm nd}$ order LDE for $\bf 
{\cal U}(N,z) $ for every $\bf N$ \\
& \\ \hline\hline
&\\
$P^{(\alpha,\beta)}_{N}$& $ \Bigl{(}z^{2}(1-z^{2})\Paeq{^{2}}{z^{2}} - 
z(2z^{2}+(\alpha -\beta)z+\alpha +\beta) \Paeq{}{z} - N(N+1+\alpha+\beta)
\Bigr{)}\ {\cal U}^{(\alpha,\beta)}(N,z) = 0 $ \\
& \\ \cline{2-2}
& \\
$T_{N}$  &  $ \Bigl{(}z^{2}(1-z^{2})\Paeq{^{2}}{z^{2}}\ +
z(1-2z^{2}) \Paeq{}{z}\ -N^{2}\Bigr{)} {\cal U}^{(T)}(N,z)\ = 0 $ \\
&\\ \cline{2-2}
&\\
$U_{N}$ & $ \Bigl{(}z^{2}(1-z^{2})\Paeq{^{2}}{z^{2}}\ -
z(1+2z^{2}) \Paeq{}{z}\ -N(N+2)\Bigr{)} {\cal U}^{(U)}(N,z)\ = 0 $\\
& \\ \hline
& \\
$Ls^{(\alpha)}_{N}$ & $ \Bigl{(}z^{3}\  \Paeq{^{2}}{z^{2}}\ + 
z(N+(1-\alpha)z)\ \Paeq{}{z}\ +
N^{2} \Bigr{)}\ {\cal U}^{(\alpha)}(N,z)\ = 0 $\\
&\\ \hline
&\\
$Hs_{N}$& $ \Bigl{(}z^{4}\  \Paeq{^{2}}{z^{2}}\ + 2z(z^{2}+N)\Paeq{}{z}\ +
2N^{2} \Bigr{)}\  {\cal U}^{(H)}(N,z)\ = 0 $\\
&\\ \hline
\end{tabular}
\end{center}
\bigskip\noindent
The above given {\it Riccati}-equations can also be found starting with 
{\it Case}'s eqs. which have been derived for the sums of powers of (simple) zeros 
of polynomials satisfying certain differential equations \cite {Ca 80}.  
{\it Table 4} shows these eqs. for the classical orthogonal polynomials. 
In the generalized {\it Laguerre} and the {\it Hermite} case scaled 
moments are used in accordance 
with the above defined scaled polynomials $Ls_{N}^{(\alpha)}(x)\ $ and 
$Hs_{N}(x)\ $:
\begin{equation}
q_{r}^{(H)}(N)\:= m_{r}^{(H)}(N)/N^{r/2}\ \ \ ,\ \ q_{r}^{(\alpha)}(N)\:=
m_{r}^{(\alpha)}(N)/N^{r}\ .  
\end{equation}
The input moments are:\par\noindent
{\it Jacobi:\ } $m_{0}(N)\equiv 1\ $\ ,\ 
$m_{1}(N):= (\beta-\alpha)/(2N+ \alpha+\beta)\ $\ .\par\smallskip\noindent
For $T_{N}:\  \alpha=\beta\ =-1/2\ $\ \ and for\ \ $U_{N}:\ \alpha=
\beta\ = +1/2\ . $\par\smallskip\noindent
generalized {\it Laguerre:\ } $q_{0}^{(\alpha)}(N)\equiv 1\ . $\par\smallskip\noindent
{\it Hermite:\ } $q_{0}^{(H)}(N)\equiv 1\ \ , \ \ q_{1}^{(H)}(N)\equiv 0\ .$
\par\smallskip\noindent
From the {\it Riccati-}eqs. which the generating functions of the power sums 
satisfy, one infers the $N\to \infty$ limits (the leading term of the $1/N$ 
expansion) which are given, together with 
the (scaled) moments, in {\it table 5}. 
$G_{0}(z):= \lim_{N\to\infty}\ G^{(\alpha,\beta)}(N,z)\ $ in the {\it Jacobi-}case.
In the generalized {\it Laguerre}-case $Q_{0}(z)\ :=\lim_{N\to\infty}\ 
Q^{(\alpha)}(N,z)\ $, and in the {\it Hermite-}case $Q_{0}^{(H)}(z):=
Q^{(H)}(N,z)\ $. $c(z)$ is the generating function of {\it Catalan}'s numbers
$\{C_{r}\}$. Universality is manifest due to the independence of the parameters
$\alpha,\beta$, resp. $\alpha$, in the {\it Jacobi-}, resp. generalized 
{\it Laguerre}-case. These asymptotic moments have also been considered in \cite{Ne-De 79} 
( although the appearance of {\it Catalan}'s numbers has not been mentioned).
\par \noindent 
\begin{center}
{\bf Table\ 4: Case's equations for moments of zeros of classical orthogonal 
polynomials }\par\smallskip\noindent
\end{center}
\begin{center}
\begin{tabular}{|c|c|}\hline
& \\
$\bf P_{N}(x)$ & {\bf Case-eqs. for $\bf{\{m_{r}(N)\}}$, resp. 
$\bf{\{q_{r}(N)\}}$, for given $\bf N$} \\
&\\ \hline\hline
& \\
$P^{(\alpha,\beta)}_{N}$ & $ (\alpha+\beta+2N-r)\ m^{(\alpha,\beta)}_{r+1}(N)
\ = (\beta-\alpha)
\ m_{r}^{(\alpha,\beta)}(N)-
r\ m^{(\alpha,\beta)}_{r-1}(N) \ + $\\ 
&\\
& $+\ N\sum_{s=0}^{r-1}\bigl{(}
m^{(\alpha,\beta)}_{r-1-s}(N)\ m^{(\alpha,\beta)}_{s}(N)\ - 
m^{(\alpha,\beta)}_{r-s}(N)\ m^{(\alpha,\beta)}_{s+1}(N)\bigr{)}$
\\
& \\ \cline{2-2}  
& \\
$T_{N}$& $ (2N-(r+1))\ m^{(T)}_{r+1}(N)\ =-r\ m^{(T)}_{r-1}(N)\ + $\\
&\\
&$+\ N \sum_{s=0}^{r-1}\bigl{(} m^{(T)}_{r-1-s}(N) \ m^{(T)}_{s}(N)\ - 
m^{(T)}_{r-s}(N)\ m^{(T)}_{s+1}(N)\bigr{)} $ \\ 
&\\ \cline{2-2}
&\\
$U_{N}$& $ (2N-(r-1))\ m^{(U)}_{r+1}(N)\ =-r\ m^{(U)}_{r-1}(N)\ +$\\  
& \\
&$+ \ N \sum_{s=0}^{r-1}\bigl{(} m^{(U)}_{r-1-s}(N) \ m^{(U)}_{s}(N)\ - 
m^{(U)}_{r-s}(N)\ m^{(U)}_{s+1}(N)\bigr{)} $\\
& \\ \hline     
& \\
$Ls^{(\alpha)}_{N}$ & $ q_{r+1}^{(\alpha)}(N)\ = {{\alpha-r}\over{N}}\ 
q_{r}^{(\alpha)}(N)\ + \sum_{s=0}^{r}\ q_{r-s}^{(\alpha)}(N)\ 
q_{s}^{(\alpha)}(N) $\\
&\\ \hline
&\\
$Hs_{N}$  & $ q_{r+2}^{(H)}(N)\ =-{{r+1}\over{2N}}\ q_{r}^{(H)}(N)\ + 
{{1}\over{2}}\sum_{s=0}^{r} q_{r-s}^{(H)}(N)\ q_{s}^{(H)}(N) $\\
& \\ \hline
\end{tabular}
\end{center}
\bigskip 
\begin{center}
{\bf Table\ 5: (Scaled) moment generating functions in the limit 
$\bf{ N\to \infty}$ }                       
\end{center}
\par\smallskip\noindent
\begin{center}
\begin{tabular}{|c|c|c|} \hline
& &\\
${\bf P_{N}(x)\ }$ & $\bf{ G_{0}(z)\ } $, {\bf resp.} $\bf {Q_{0}(z)\ }$& 
{\bf limit moments} \\
& &\\ \hline\hline
& &\\
$P^{(\alpha,\beta)}_{N}$ &  $G_{0}(z)\ = 1/\sqrt{1-z^2}\ $& $m_{0,2l}\ = {{2l}
\choose {l}}/2^{2l}\ $ \\
& &\\ \hline
& &\\
$Ls^{(\alpha)}_{N}$ & $Q_{0}(z)\ = {{1}\over{2z}}(1-\sqrt{1-4z})\ =: c(z)\ $&
$q_{0,r}\ = C_{r}\ := {{2r}\choose{r}}/(r+1)$ \\
& &\\ \hline
& & \\
$Hs_{N}$  & $Q_{0}^{(H)}(z)\ = c(z^{2}/2)\ = {{1}\over{z^{2}}}
(1-\sqrt{1-2z^{2}})\ $&
$q_{0,2l}^{(H)}\ = C_{l}/2^{l}$ \\
& &\\ \hline
\end{tabular}
\end{center}
It is clear from the explicit expressions for the generating functions 
$G(N,z)$ in both {\it Chebyshev} instances that a $1/N$ expansion is doomed
to fail. Such an expansion produces for the $T_{N}-$polynomials the asymptotic
moment generating function $G_{0}(z)$ given for the {\it Jacobi}-case in 
{\it table} 5. This happens to be a solution of the $T_{N}-${\it Riccati}-eq. 
shown in {\it table 2}
for all $N$. However, it is not the correct solution $G^{(T)}(N,z)$ given in
{\it table 1}. This explains the mistake which occured in \cite{Ca 80}. It is 
possible to recover $G^{(T)}(N,z)$ from the {\it special} solution $G_{0}(z)$ 
of the $T_{N}-${\it Riccati}-equation. In the {\it Chebyshev} $U_{N}$ case the sum of 
the first two $1/N$-expansion coefficients furnishes a {\it special} solution 
of the $U_{N}-${Riccati-}eq. (see \cite {La 97}), which is also not the correct 
solution $G^{(U)}(N,z)\ $ given in {\it table 1}.
\par\bigskip\noindent
\section{\bf Distribution of zeros  and Perron-Stieltjes inversion }\par\smallskip
\def\Chi{\raise2.5pt\hbox{$\chi$}}
\hskip 1cm For polynomials with only real zeros the (ordinary) generating function  $G(N,z)$ 

of the sums of powers of zeros (or moments) is related to the {\it Stieltjes-}transform
$\Chi(N,z)$ of the discrete real measure $(4)$in a simple way. The {\it Stieltjes}-transform
of $\rho_{N}(x)$ is defined by
\begin{equation}
\Chi(N,z)\ := \int_{-\infty}^{+\infty} {{1}\over{z-x}}\ d\rho_{N}(x)\ \ \ ,
\ \ \ z\not\in supp(d\rho_{N})\ .
\end{equation}
Because (in the sense of formal power series) one has with the counting 
measure (4)  
\begin{equation}
G(N,z)\ =  \sum_{r=0}^{\infty} m_{r}(N)\  z^{r} \ = 
{{1}\over{N}}\sum_{i=1}^{N}{{1}\over{1-z\ x^{(N)}_{i}}} \ =
\int_{-\infty}^{\infty}{{1}\over{1-z\ x}}\ d\rho_{N}(x) \ , 
\end{equation}
one finds the simple relation 
\begin{equation}
\Chi(N,z)\ = {{1}\over{z}}\ G(N,{{1}\over{z}})\ \ .
\end{equation}
It is possible to recover the discrete measure from $\Chi(N,z)$, 
{\it i.e} $G(N,z)$, with the aid of the {\it Perron-Stieltjes} inversion 
formula \cite{Wi 29},\cite{Ak 65},\cite{Ch 78},\cite{As 84} 
\begin{equation}
\rho_{N}(t_{2})\ -\  \rho_{N}(t_{1})\ =-{{1}\over{\pi}}\ \lim_{\eta\to +0}\ 
\int_{t_{1}}^{t_{2}}dt\ Im\ \Chi(N,t+i\eta)\ \ ,
\end{equation}
with $\rho_{N}(t_{k})\ :={{1}\over{2}}(\rho_{N}(t_{k}+0)\ +\rho_{N}(t_{k}-0))$
for $k=1,2$. This inversion is valid for real measures of bounded variation,
provided ${\bar {\Chi}}_{N}=\Chi_{N}({\bar z})$ for $z\in {\bf C}$, and
$\Chi_{N}$ is analytic off the real axis. As long as $N$ is finite, 
$\Chi_{N}(z)$ has 
poles only at the zeros of $P_{N}$. The inversion formula then yields for 
$\rho_{N}$ jumps at these real zeros $x_{i}^{(N)}$ with height 
$m_{i}^{(N)}/N$ where $m_{i}$ counts the multiplicity of the zero $x_{i}^{(N)}$. 
($(36)$ shows that poles in  $\Chi$ lead to jumps in $\rho$. 
The converse statement is obvious from $(33)$.) In this way one recovers 
the purely discontinuous measure $(4)$.\par\smallskip\noindent
The inversion formula becomes especially interesting if the limit 
$\Chi(z)\ :=\lim_{N\to \infty} \Chi(N,z)$ exists. Then one can find the limit 
measure or distribution of zeros. In our examples of classical (scaled) 
orthogonal polynomials the limit moment generating functions exist 
(see {\it table 5}), and hence one can find by {\it Perron-Stieltjes} inversion
of $\Chi(z):= {{1}\over{z}}\ G_{0}({{1}\over{z}})\ $ (resp. ${{1}\over{z}}\ 
Q_{0}({{1}\over{z}}) $) the distribution of zeros in the limit $N\to \infty$. 
These are
listed in {\it table 6}. In all cases the conventional definition of the 
square-root in $\bf C$ has been used. These limit distributions (with other 
normalization) have been considered in \cite{GaDe 86}.
\begin{center}
{\bf Table\ 6: Limit distributions of zeros of (scaled) orthogonal 
polynomials }                       
\end{center}
\par\smallskip\noindent
\begin{center}
\begin{tabular}{|c|c|} \hline
& \\
${\bf P_{N}(x)\ }$ & $\rho(x)\ =\ \lim_{N\to\infty} \rho_{N}(x)$ \\
& \\ \hline\hline
& \\
$P^{(\alpha,\beta)}_{N}$ &  $\rho^{(\alpha,\beta)}(x)\ = {{1}\over{\pi}}\ 
{{1}\over{\sqrt{1-x^{2}}}}\ \ , \ \ {\rm for}\ \vert x \vert \ < \ 1\ {\rm and} \ 0\ 
{\rm otherwise}\ $ \\
& \\ \hline
& \\
$Ls_{N}^{(\alpha)}$&$ \rho^{(Ls^{(\alpha)})}(x)\ = {{1}\over{2\pi}} \sqrt{{{4}\over{x}}-1}\ ,\  
{\rm for} \  0< x \leq 4\ \ , \ \  {\rm and\ } 0 {\ \rm otherwise}$ \\
& \\ \hline
& \\
$Hs_{N}$  & $\rho^{(Hs)}(x)\ = {{1}\over{\pi}} \sqrt{2-x^{2}}\ \ ,\ \  {\rm for}\ \
\vert x\vert <\sqrt{2}\ \ {\rm and\ } 0 {\ \rm otherwise}$\\ 
& \\ \hline
\end{tabular}
\end{center}
\section{Conclusion}
\hskip 1cm Formula $(25)$ which gives the (ordinary) generating function 
$G(N,z)$ for the sums of (positive) powers of zeros of any one-variable 
polynomial $(1)$ can be simplified whenever the (logarithmic) derivative 
rule for the polynomial under consideration is known. If the polynomial 
satisfies a (second) order differential equation (like in the assumption
of the work  \cite{Ca 80}) one expects $G(N,z)$ to satisfy a 
{\it Riccati-}type differential equation. For polynomials with only real 
zeros ({\it e.g.} orthogonal polynomials with positive moment functionals)
one can use {\it Perron-Stieltjes} inversion to find the distribution of 
zeros in the limit $N\to \infty$, provided $\lim_{N\to\infty}\ G(N,1/z)/z\ $
exists.\par\bigskip\noindent
\section *{ Acknowledgements}
\hskip 1cm The author likes to thank Professor J. Dehesa for acquainting him 
with the work of the Granada group on sums of powers of zeros of polynomials
by sending him his reprints. He also would like to thank the organizers of the
Sevilla Symposium on Orthogonal Polynomials for the stimulating atmosphere.

\vfill\eject
\end{document}